\documentclass[11pt,twoside]{article}
\usepackage{latexsym,amsmath}
\usepackage{indentfirst}
\usepackage{graphicx}

\begin{document}
\def \b{\Box}
\begin{center}
{\Large {\bf Dynamical behaviors of the special fractional-order Chen-Lee system}}
\end{center}

\begin{center}
{\bf Mihai Ivan}\\
\end{center}

\setcounter{page}{1}
\pagestyle{myheadings}

{\small {\bf Abstract}. The main purpose of  this paper is to study the special fractional-order Chen-Lee system, using the Caputo fractional derivatives. For this fractional model we investigate the existence and uniqueness of solution of initial value problem, 
asymptotic stability of its equilibrium states,  stabilization problem using appropriate control and numerical integration.
{}

\large {\bf 1. Introduction}\\

{\small The theory of fractional differential equations (i.e. fractional calculus) and its applications are based on  non-integer order of derivatives and integrals \cite{kilb, podl,  suzb}.

Fractional calculus has found many applications in various areas of science and engineering. In the last decades, studies were carried out in research fields leading to the exploration of the dynamic behaviors of some classical differential systems \cite{enpp,  igmp, labi, lima, pagi}, as well as some fractional-order differential systems. For example, the fractional models played an important role in  applied mathematics \cite{bora, gi22,  iomi, mlag, migo}, mathematical physics \cite{nabu, igom, miads, miva23}, applied physics \cite{hilf, mist, miisj, sczc},  mathematical biology \cite{acbv, ahme, iyoa, lizh, kukc}, control processing and nonlinear circuits \cite{mati, podl, wlccc}, control chaos and synchronization \cite{ivmi, mato} and so on.

The Chen-Lee system \cite{chle} is described by the following differential equations on $ {\bf R}^{3}: $\\[-0.2cm]
\begin{equation}
 \dot{x}^{1}(t)  =
 - x^{2}(t) x^{3}(t) + a  x^{1}(t),~~ \dot{x}^{2}(t)  =
  x^{1}(t)x^{3}(t) - b  x^{2}(t) ,~~ \dot{x}^{3}(t) = \frac{1}{3} x^{1}(t)x^{2}(t) - c  x^{3}(t) ,\label{(1)}
\end{equation}
 where  $~x^{i}~$ are state variables, $ \dot{x}^{i}(t) = d x^{i}(t)/dt~$ for $ i=\overline{1,3}~$ and $~a,b,c\in {\bf R}~$ are parameters.

Behavior of many dynamical systems can be described and studied using the associated fractional models. In this context, we associate to system (1), a fractional model, denoted by (11). This system is called the  fractional-order Chen-Lee system and it is obtained from the system  $ (1), $  replacing the ordinary differential operators by the fractional-order differential ones.

The fractional system $ (11) $  with Caputo's fractional derivatives of order $ q (q\in (0, 1]) $  has proven to be a useful tool in synchronization of two chaotic systems and implementing electronic circuits. For this purpose, the fractional systems  $ (11), $ determined by non-zero real values given to parameters $ a, b, c, $  have discussed in \cite{shcct, stcl, wlccc}.

In the paper \cite{ivan24} was presented  some applications of groupoids  in the theory of classical and fractional dynamical systems, 
 using  Lie groupoids, Leibniz algebroids and Lie algebroids (\cite{comgimg09, imod, migi, migopr, migo}).

In this paper we define a new fractional system, denoted by $ (12) $ and  called the {\it special fractional-order Chen-Lee system}. The dynamics of this fractional model is described by equations $ (11) $ in which the second parameter has been replaced by the zero value.

The structure of this paper is as follows. In Section 2 we introduce some notations, definitions and preliminary facts about fractional differential systems which are useful throughout the paper. In Section 3, the existence and uniqueness of solution of initial value problem for the special fractional Chen-Lee system  is proven. Also,  the asymptotic stability of equilibrium states for this fractional model is analysed. Section $ 4 $ is devoted to establishing sufficient conditions that the parameters $ k, a $ and $ c $ satisfy to control the chaotic behavior of the special fractional Chen-Lee system with control around its equilibrium states. In Section $ 5,$ we present the numerical integration of fractional model with control  $(15),$ using the first fractional Euler's scheme. Finally, a numerical simulation of $~(15)~$ is presented to illustrate the theoretical results.}\\

\large {\bf 2. Preliminaries on fractional differential systems }\\

We present some definitions and basic results related to Caputo fractional derivatives that are needed in this paper.

For the qualitative study of fractional differential systems, the derivative of fractional-order $ q>0 $ in the sense of Caputo is often used. In many concrete applications, the study of fractional systems is based on the properties of Caputo's \emph{}derivative of fractional-order  $ q\in (0,1). $
In this paper we suppose that $ q \in (0,1].$

 Let $ f\in C^{\infty}(\textbf{R}) $ and $ \alpha \in (0, +\infty).~$ The $ \alpha-$order Riemann-Liouville integral operator, denoted by $~I^{\alpha}f(t),~$ is defined as:\\[-0.2cm]
 \begin{equation}
I^{\alpha}(t) = \displaystyle\frac{1}{\Gamma(\alpha)}\int_{0}^{t}{(t-s)^{\alpha
-1}}f(s)ds,~\alpha> 0,\label{(1)}
\end{equation}
  where $~\Gamma~$ is the Euler's Gamma  function.

The Caputo differential operator of order $~ q> 0~$ (\cite{difo}), is defined as:\\[-0.4cm]
\begin{equation}
D_{t}^{q}f(t) = I^{m-q}f^{(m)}(t),~~~q > 0.\label{(2)}
\end{equation}
where $~f^{(m)}(t)$  represents the $ m-$order derivative of the function $ f~$ and $~m \in \textbf {N}^{\ast}$ is an integer such that $ m-1 \leq q \leq m.~$  If $ q=1,$  then $ D_{t}^{q}f(t) = \frac{df}{dt}.$
\markboth{\it M. Ivan}{ Dynamical behaviors of the special fractional-order Chen-Lee system}

According to relations $ (2) $ and $ (3),$ one obtains the following definition.

{\bf Definition 1.} {\rm The Caputo fractional derivative of order $~q >0~$ of the function $ f\in C^{\infty}({\bf R}) $ is given by:\\[-0.4cm]
\begin{equation}
D_{t}^{q}f(t) = \displaystyle\frac{1}{\Gamma(m-q)}\int_{0}^{t}{(t-s)^{m-q-1}}f^{(m)}(s)ds,~~~q > 0.\label{(3)}
\end{equation}

In particular, the Caputo derivative of order $ q\in (0,1) $ of $ f\in C^{\infty}({\bf R}) $ is described by:\\[-0.2cm]
\begin{equation}
D_{t}^{q}f(t) = \displaystyle\frac{1}{\Gamma(1-q)}\int_{0}^{t}{(t-s)^{-q}}f'(s)ds,~~~q \in (0, 1).\label{(4)}
\end{equation}

 In the Euclidean space $ {\bf R}^{n} $ with local coordinates $~\{ x^{1}, x^{2}, \ldots, x^{n}\}, $ we consider the following system of fractional differential equations:\\[-0.2cm]
\begin{equation}
D_{t}^{q}x^{i}(t) = f_{i}(x^{1}(t), x^{2}(t), \ldots,
x^{n}(t)) ,~~ i=\overline{1,n},\label{(5)}
\end{equation}
where $q\in (0,1), f_{i}\in C^{\infty}({\bf R}^{n}, {\bf R}),
~ D_{t}^{q} x^{i}(t)$ is the Caputo derivative of order $ q$ for $i=\overline{1,n}$ and $t\in [0,\tau)$ is the
time.

The fractional system $(6)$ can be written as follows:
\begin{equation}
D_{t}^{q}x(t) = f(x(t)) ,\label{(6)}
\end{equation}
where $~f(x(t)) = (f_{1}(x^{1}(t),\ldots, x^{n}(t)),
f_{2}(x^{1}(t),\ldots, x^{n}(t)), \ldots, f_{n}(x^{1}(t),\ldots,
x^{n}(t)))^{T} $ and $ D_{t}^{q} x(t)= ( D_{t}^{q}
x^{1}(t), \ldots, D_{t}^{q} x^{n}(t))^{T}.$

Consider the initial value problem with Caputo derivative for the system $ (7) $ :\\[-0.2cm]
\begin{equation}
D_{t}^{q} x(t) = f(x(t)),~~~ x(0)=x_{0},~~~t\in I=[0,T],~T>0  \label{(8)}
\end{equation}
where $~x: I \rightarrow {\bf R}^{n},~f: {\bf R}^{n} \rightarrow {\bf R}^{n}~$ is a continuous function and
$ q\in (0,1).$

{\bf Theorem 1.} {\rm (\cite{acbv, difo})} {\it  Let be the initial value problem given by  Equation $ (8).$  If the function $ f $ satisfies the following two conditions:\\
$(i)~f~$ is differentiable and bounded  on $~ D \subset {\bf R}^{n};$\\
$(ii)~f(x(t))~$ satisfies the local Lipschitz condition with respect to $ x $,\\
then there exists a unique solution of Equation $ (8) $ on the $~[0, T]\times D. ~~~\hfill\Box $}

A point $x_{e}=(x_{e}^{1}, x_{e}^{2},\ldots, x_{e}^{n})\in {\bf
R}^{n}$ is said to be {\it equilibriun point} or {\it equilibrium state} of system
$(6)$, if $~D_{t}^{q}x^{i}(t) =0$ for $i=\overline{1,n}$.

To analyse the nature of equilibrium points of the system $ (6), $  the Jacobian matrix $ J(x) $ is determined. The matrix $ J(x), $  associated with  $(6),$ is given by:\\[-0.2cm]
\begin{equation}
J(x)=(\displaystyle\frac{\partial f_{i}}{\partial
x^{j}}),~~~~~i,j=\overline{1,n}. \label{(9)}
\end{equation}

Denote with $~\lambda_{i},~i=\overline{1,n},~$ the eigenvalues of matrix $~J (x_{e}),$ where $~J (x_{e})$ is the Jacobian matrix $J(x)$ evaluated at  equilibrium point $ x_{e}. $

{\bf Proposition 1.} {\rm (\cite{mati}) ({\bf Matignon's test}) {\it Let $ x_{e} $ be an equilibrium point  of fractional system $(6)$ and
$~\lambda_{i}, i=\overline(1,n),~ $ the eigenvalues of matrix $~J(x_{e}).$

 $(i)~ x_{e}$ asymptotically stable, if and only if  the following relations are satisfied:\\[-0.2cm]
\begin{equation}
| arg(\lambda_{i}) | > \displaystyle\frac{q\pi}{2}, ~~~i=\overline{1,n}.\label{(10)}
\end{equation}
 $(ii)~ x_{e} $ is locally stable, if and only if either it is asymptotically stable, or the
eigenvalues $ \lambda_{j} $ satisfying $~| arg(\lambda_{j}) | = \displaystyle\frac{q \pi}{2}~$ have geometric
multiplicity one.}\hfill$\Box$

{\bf Lemma  1.} {\it Let $~x_{e}~$ be an equilibrium state of the fractional model $~(6),~\lambda_{i},~i=\overline{1,n}~$ the eigenvalues of $~J(x_{e})~$
and $ q\in (0,1).$

$(i)~$ If one of the eigenvalues $~\lambda_{i},~i=\overline{1,n}~$ is equal to zero or it is positive, then $~x_{e}~$ is unstable.

$(ii)~$ If $~\lambda_{i} < 0, $ for all $~i=\overline{1,n},~$ then  $~x_{e}~$   is asymptotically stable}.

{\bf Proof.} $~(i) $  In the case when there exists one eigenvalue $~\lambda_{j},~j\in \{1,2, \ldots, n\} $ such that $~\lambda_{j}=0 ~$ or $~ \lambda_{j} > 0, $ then $~| arg(\lambda_{j}) | =0 < \displaystyle\frac{q \pi}{2}, (\forall) q\in (0,1).~$ Acording to Proposition 1(i), it follows that $~x_{e} $ is unstable.

$~(ii) $ Assume $~\lambda_{i}  < 0, i=\overline{1,n}.$ Then $~| arg(\lambda_{i}) | = \pi > \displaystyle\frac{q \pi}{2}, (\forall) q\in (0,1).~$ By Proposition 1(i), it follows that $~x_{e} $ is assymptotically stable.\hfill$\Box$\\

\large {\bf 3. The special fractional-order Chen-Lee system}\\[-0.4cm]

{\bf 3.1. Description of the special fractional model}

In this subsection we define the special fractional-order Chen-Lee system. For this new fractional model,  the existence and uniqueness of solution of initial value problem are proved.

The Chen-Lee system (\cite{chle}) was derived from the Euler equations for the motion of a rigid body with principal axes at the center of mass.

The fractional model associated to Chen-Lee system $ (1), $ called the {\it fractional Chen-Lee system}, is defined by the following fractional differential equations:\\[-0.2cm]
\begin{equation}
\left\{ \begin{array} {lcl}
 D_{t}^{q}{x}^{1}(t) & = &
 - x^{2}(t) x^{3}(t) + a x^{1}(t)  \\[0.1cm]
 D_{t}^{q}{x}^{2}(t) & = & x^{1}(t) x^{3}(t)- b  x^{2}(t),~~~~~a, b, c\in{\bf R},~~~~~  q \in (0,1),\\[0.1cm]
  D_{t}^{q}{x}^{3}(t) & = & \frac{1}{3} x^{1}(t)x^{2}(t) - c x^{3}(t). \label{(11)}
  \end{array}\right.
\end{equation}

We will associate to Equations $ (11), $  the parameter vector $ v_{CL} = (a,b,c). $  If $ a, b $ and $ c $ are replaced with arbitrary real values, then the vector $ v_{CL} $ generates a family of fractional systems of the  Chen-Lee type.

{\bf Remark 1.} {\rm  The fractional systems of the Chen-Lee type for $ q=1 $ or $ q\in (0,1), $ have been studied from various research directions by many authors.

{\bf 1.} From the point of view of Poisson geometry, the Chen-Lee system $ (11) $  for $~q=1, $ was researched in (\cite{labi, popa}). More preciselly:\\
$(i)~$ in \cite{popa} are studied  the systems $ (11) $ determined by $ q=1 $ and the vectors $~v_{CL}^{1} = (a, 0, 0), v_{CL}^{2} = (0, b, 0), v_{CL}^{3} = (0, 0, c) $ with $ a,b,c\in {\bf R}^{\ast}.$\\
$(ii)~$ in \cite{labi} are analyzed the symmetries of a class of systems on $ {\bf R}^{3}. $ The Chen-Lee system $ (11) $ with $~ q=1 $ and  $~v_{CL}^{0} = (0, 0, 0) $   belongs of considered class.

{\bf 2.} Another direction of research is oriented towards studying the chaotic behavior of classical or fractional dynamical systems. The chaotic behavior of fractional systems of Chen-Lee type has investigated by changing the initial conditions and the associated vector $ v_{CL}, $  using different methods and numerical simulations. For instance:\\
$(i)~$ in \cite{shcct} is studied the synchronization of two Chen-Lee systems determined by $ q=1 $ and the constant  vector $~v_{CL}^{4} = (5, -10, -3.8). $\\
$(ii)~$ in \cite{stcl}, two Chen-Lee systems,  determined by $~q=1,~v_{CL}^{4} $ and $~v_{CL}^{5} = (3, -5, -1),~$ are realized as nonlinear circuits.\\
$(iii)~$ in \cite{wlccc} is investigated the fractional Chen-Lee system $ (11) $  with the vector $~v_{CL}^{4}, $ for $ q=0.9. $  Also, this chaotic system is realized  as a electronic circuit.}\hfill$\Box$

If in $ (11) $ we take $ ~v_{CL} = (a, 0, c), $ one obtains a new fractional system. This new system is defined by the following fractional differential equations:\\[-0.2cm]
\begin{equation}
\left\{ \begin{array} {lcl}
 D_{t}^{q}{x}^{1}(t) & = &
 - x^{2}(t) x^{3}(t) + a x^{1}(t)  \\[0.1cm]
 D_{t}^{q}{x}^{2}(t) & = & x^{1}(t) x^{3}(t)  ,~~~~~~~~~~~~~~~~  q \in (0,1),\\[0.1cm]
  D_{t}^{q}{x}^{3}(t) & = & \frac{1}{3} x^{1}(t)x^{2}(t) - c x^{3}(t), \label{(12)}
  \end{array}\right.
\end{equation}
 where $~a, c\in{\bf R}~$ such that $ ac \neq 0.$

The system $ (12) $ is called the {\it special fractional Chen-Lee system}. It is determined by the vector $ v_{CL} $ which has two non-zero components.

The initial value problem for the special fractional Chen-Lee system
$ (12) $ can be represented in the following matrix form:\\[-0.2cm]
\begin{equation}
D_{t}^{\alpha}x(t)  =  x^{1}(t) A x(t) +  B  x(t) ,~~~~~ x(0) =
x_{0},\label{(13)}
\end{equation}
where $0 < q < 1,~ x(t)= ( x^{1}(t),
 x^{2}(t), x^{3}(t))^{T}, ~t\in(0,\tau)$ and\\[-0.1cm]
\[
A = \left ( \begin{array}{ccc}
0 & 0 & 0 \\
0 & 0 & 1 \\
0 & \frac{1}{3} & 0 \\
\end{array}\right ),~~~ B = \left ( \begin{array}{ccc}
a & 0 & 0 \\
1 & 0 & 0\\
0 & 0 & -c\\
\end{array}\right ).
\]
{\bf Theorem  2.} {\it The initial value problem for the special fractional Chen-Lee system
$ (12) $ has a unique solution}.

{\bf Proof.} Let $ f(x(t))= x^{1}(t) A x(t) + B x(t).$ It is
obviously continuous and bounded on $ D =\{ x \in {\bf R}^{3} |~
x^{i}\in [x_{0}^{i} - \delta, x_{0}^{i} + \delta]\}, i=\overline{1,3} $ for any
$\delta>0. $ We have $~f(x(t)) - f(y(t)) =  x^{1}(t) A x(t) - y^{1}(t) A y(t) + B x(t)-  B y(t)=g(t)+h(t), $
where $~g(t)= x^{1}(t) A x(t) - y^{1}(t) A y(t)~$ and $~h(t)= B x(t) - B y(t). $
 Then\\[0.2cm]
$(a)~~|f(x(t)) - f(y(t))|\leq |g(t)| + |h(t)|. $

It is easy to see that $~g(t)= (x^{1}(t)- y^{1}(t))A
x(t)+ y^{1}(t) A(x(t)- y(t)).~$ Then\\[0.2cm]
$|g(t)| \leq |(x^{1}(t)- y^{1}(t)) A x(t)| + |y^{1}(t)
A (x(t)- y(t))|\leq $\\
$\leq  \|A\|( |x(t)|\cdot |x^{1}(t)- y^{1}(t)| + |y^{1}(t)|\cdot |x(t) -y(t)|, $\\
where $ \|\cdot \| $ and $|\cdot|$ denote matrix norm
and vector norm respectively.

Using the inequality $~|x^{1}(t)- y^{1}(t))|\leq |x(t)- y(t))|~$ one obtains\\[0.1cm]
$(b)~~~|g(t)| \leq (\|A\|+ |y^{1}(t)|)\cdot |x(t)- y(t)|.$\\[0.1cm]
Similarly, we have\\[0.2cm]
$(c)~~~|h(t)| \leq \|B\|\cdot |x(t)- y(t)|.$\\[0.1cm]
According to $(b)$ and $ (c),$ the relation $(a)$ becomes\\[0.1cm]
$(d)~~~|f(x(t)) - f(y(t))|\leq  (\|A\| + \|B\| +
|y^{1}(t)|)\cdot |x(t)-y(t)|.$\\[0.1cm]
Replacing  $\|A\|= \frac{\sqrt{10}}{3},~ \|B\| = \sqrt{1+c^{2}}$ and using the inequalities $~|y^{1}(t)|\leq |x_{0}| + \delta,~$ from the
relation $ (d), $ we deduce that\\[0.1cm]
$(e)~~~|f(x(t)) - f(y(t))|\leq  L\cdot |x(t)-y(t)|,~~~
\hbox{where}~L =1 + \frac{\sqrt{10}}{3}+\sqrt{1+c^{2}} + |x_{0}| + \delta > 0.$

 The inequality $(e)$ shows that $ f(x(t))$
satisfies a Lipschitz condition. Based on the results of Theorem
$ 1, $  we can conclude that for the initial value problem of the system $ (12) $ there exists  a unique solution. \hfill$\Box$\\[-0.2cm]

{\bf 3.2. Stability analysis  of the special fractional model}

In this subsection we present the study  of asymptotic stability for the equilibrium states of  fractional model $ (12). $
 For this study we apply the Matignon's test.

For the system $ (12) $ we introduce the following
 notations:\\[-0.2cm]
 \begin{equation}
f_{1}(x) =  - x^{2}x^{3} + a  x^{1},~~~ f_{2}(x) =
x^{1} x^{3},~~~ f_{3}(x)  = \frac{1}{3} x^{1} x^{2} -c  x^{3}.\label{(14)}
\end{equation}
{\bf Proposition 2.} {\it The equilibrium states of the special fractional-order Chen-Lee system
$(12)$ are given as the following
family}:\\[-0.3cm]
\[
E_{2}:=\{ e_{2}^{m}=(0, m, 0)\in {\bf R}^{3} |~ m \in {\bf
R}\}.
\]
{\bf Proof.} The equilibrium states are solutions of the equations
$~f_{i}(x)=0, i=\overline{1,3}$ where $~f_{i},~i=\overline{1,3}$
are given by (14).\hfill$\Box$

The Jacobian matrix associated to system $ (12) $ is:
\[
J(x) = \left ( \begin{array}{ccc}
a & -x^{3}    & -x^{2} \\
  x^{3}   & 0 &  x^{1}\\
  \frac{1}{3} x^{2}   &  \frac{1}{3} x^{1}  & -c \\
\end{array}\right ).\\[-0.1cm]
\]

{\bf Proposition 3.} {\it The equilibrium states $ e_{2}^{m}\in E_{2}~$
are unstable $ (\forall) q \in (0,1).$}

{\bf Proof.} The characteristic polynomial of $~
J(e_{2}^{m}) =\left (\begin{array}{ccc}
  a & 0  & -m\\[0.1cm]
  0 & 0 & 0 \\
  \frac{m}{3} & 0  & -c \\
\end{array}\right )~ $
is\\
 $~ p_{J(e_{2}^{m})}(\lambda) = \det ( J(e_{2}^{m}) -
\lambda I) = - \lambda [ (\lambda -a)(\lambda +c) + \frac{1}{3} m^{2}].~$
 The equation $ ~p_{J(e_{2}^{m})}(\lambda) = 0 $  has the root $ \lambda_{1} = 0 .$  Since $~arg(\lambda_{1}) =0 < \frac{q \pi}{2}$
for all $ q\in (0,1),$  by Lemma 1(i), follows that $ e_{2}^{m}, m\in {\bf R}~$ are unstable for all $ q\in (0,1).$ \hfill$\Box$\\[-0.2cm]

\large {\bf 4.~ Controlling the chaotic behavior  of the special fractional model }\\[-0.2cm]

In order to control the chaotic behavior of the special fractional model $ (12), $ we will build a fractional system with control associated with this system.\\[-0.3cm]

{\bf 4.1.  Stability analysis of the special fractional Chen-Lee system with control}\\[-0.3cm]

To stabilize the chaotic behavior of a fractional model, we will apply the method of stabilizing unstable equilibrium states of a given system (\cite{giva}).

Let $~x_{e}= (x_{e}^{1}, x_{e}^{2}, x_{e}^{3}) $  be an unstable equilibrium state of the fractional system $(12). $  We associate with  $(12)$ a new
fractional system, called the {\it special fractional Chen-Lee system with control} at  $ x_{e}. $  This system is given by:\\[-0.2cm]
\begin{equation}
\left\{ \begin{array} {lcl}
 D_{t}^{q}{x}^{1}(t) & = &
 - x^{2}(t) x^{3}(t) + a x^{1}(t) \\[0.1cm]
 D_{t}^{q}{x}^{2}(t) & = & x^{1}(t) x^{3}(t) + k (x^{2}(t)- x_{e}^{2}) ,~~~~~~~~~~~~~~~~~ q \in (0,1),\\[0.1cm]
  D_{t}^{q}{x}^{3}(t) & = & \frac{1}{3} x^{1}(t) x^{2}(t) - c x^{3}(t) \label{(15)}
  \end{array}\right.
\end{equation}
where  $~a, c\in{\bf R}^{\ast}~$ are real constants and $~k\in {\bf R}^{\ast}~$ is a control parameter.

The Jacobian matrix of the fractional model  with control $(15),$ is given as:\\[-0.2cm]
\[
J_{k}(x, a, c) = \left ( \begin{array}{ccccc}
a & - x^{3} & - x^{2} \\
  x^{3}  & k &   x^{1}\\
   \frac{1}{3} x^{2}   &  \frac{1}{3} x^{1}  & - c \\
\end{array}\right ).
\]

If one selects the parameters $ a, c, k $ which then make the eigenvalues of the Jacobian matrix of fractional model
$(15)$ satisfy the condition from Proposition 1, then its trajectories asymptotically approaches the unstable
equilibrium state $ x_{e} $ in the sense that $\lim_{t\rightarrow \infty} \|x(t)- x_{e}\|= 0$, where $\|\cdot\|$ is the Euclidean norm.

{\bf Theorem 3.} {\it Let  be the special fractional Chen-Lee system with the control $~k \in {\bf R^{\ast}},~ q\in (0,1)~$ and $~e_{0}= (0,0,0).$\\
$~~~~~~~$ {\bf 1.} $~k < 0.$\\
$(i)~$ If $~ a< 0 $  and $~ c >0, $ then $ e_{0}~$ is asymptotically stable.\\
$(ii)~$ If $~a c > 0~$  or $~ a >0 $  and $~ c <0,~ $ then $~ e_{0}~$ is unstable.\\
$~~~~~~~$ {\bf 2.} $~k > 0.~$ If  $~a, c\in {\bf R}^{\ast},~$ then $~ e_{0} $ is unstable.}

{\bf Proof.} The characteristic polynomial of $~J_{k}(e_{0}, a, c) $ is  denoted with $~p_{0}(\lambda),~$ where
 $~p_{0}(\lambda) =-(\lambda - a)(\lambda + c)(\lambda -k).$ The roots of equation $~p_{0}(\lambda)= 0~$ are\\
 $\lambda_{1}= a,~\lambda_{2}=-c,~\lambda_{3}= k.$\\
{\bf 1.} {\bf Case} $~k <0.~$ Then $~\lambda_{3} <0.~$\\
$(i)~$ We have $~\lambda_{1}<0~$ and $~\lambda_{2}<0~$  if and only if $~\lambda_{1} + \lambda_{2} <0 $ and  $~\lambda_{1}\cdot \lambda_{2} >0.$ It follows that $~\lambda_{i} <0, i=\overline{1,3}~$ for all $ a,c\in {\bf R}^{\ast}~$ such that $~a<0 $ and $ c >0.$  According to Lemma 1(ii),
$ e_{0}~$ is asymptotically stable  for all $~k\in (-\infty, 0).$\\
$(ii)~$ We suppose $~a c > 0~$  or $~ a >0 $  and $~ c <0.~ $ In this case, one from the eigenvalues $~\lambda_{1}~$  and $~\lambda_{2}~$ is positive. Since $~J_{k}(e_{0},a, c)$ has at least a positive eigenvalue, by Lemma 1(i), it follows that  $ e_{0} $ is unstable.\\
{\bf 2.} {\bf Case} $~k >0.~$ Since $~J_{k}(e_{0},a, c)$ has at least a positive eigenvalue, it follows that  $ e_{0} $ is unstable. \hfill$\Box$

{\bf Theorem 4.} {\it Let  be the  fractional Chen-Lee system with the control $~k \in {\bf R^{\ast}}. $ Let $~\Delta = (a+c)^{2} - \frac{4}{3} m^{2},~ q_{2} = \displaystyle\frac{2}{\pi}
\arctan\displaystyle\frac{\sqrt{-\Delta}}{a -c}~$ and $~e_{2}^{m}= (0,m,0),~m\neq 0.$\\
$~~~~~~~$ {\bf 1.} Let $~k < 0~$ and $~\Delta < 0.$\\
$(i)~$ If  $~a<c,~$ then $~e_{2}^{m}$ is asymptotically stable  $~(\forall)~ m \in (-\infty, -\frac {|a+c|\sqrt {3}}{2})\cup ( \frac {|a+c|\sqrt {3}}{2}, \infty )~$ and $~ q\in (0,1).$\\
$(ii)~$ If  $~a  > c~$ and $~m \in (-\infty, -\frac {|a+c|\sqrt {3}}{2})\cup ( \frac {|a+c|\sqrt {3}}{2}, \infty ),~$ then:\\
{\bf (a)}$~~ e_{2}^{m}~$ is asymptotically stable  $~(\forall) q \in ( 0, q_{2})~$ and it is stable for $~q = q_{2};$\\
{\bf (b)} $~~ e_{2}^{m}~$ is unstable $~(\forall)~q \in ( q_{2}, 1).$\\
$~~~~~~~$ {\bf 2.} Let $~k < 0~$ and $~\Delta > 0.$\\
$(i)~$ If  $~a c > 0~$ and $~a < c,~$ then $~e_{2}^{m}~$ is asymptotically stable  for
all  $~m\in (-\frac {|a+c|\sqrt {3}}{2},-\sqrt{3ac})\cup (\sqrt{3ac}, \frac {|a+c|\sqrt{3}}{2})~$ and $~ q\in (0,1).$\\
$(ii)~$ If  $~a c < 0,~$ then $ e_{2}^{m}~$ is unstable  for all $~m \in (-\sqrt{3ac}, \sqrt{3ac})~$ and $~ q\in (0,1).~$\\
$~~~~~~~$ {\bf 3.} Let $~k > 0.~$ If  $~a, c\in {\bf R}^{\ast},~$  then $~ e_{2}^{m}$ is unstable $~(\forall) m\in {\bf R}^{\ast}~$ and $~q \in (0, 1).$}

{\bf Proof.} The characteristic polynomial of matrix $~J_{k}(e_{2}^{m}, a, c),~$ is \\
$p_{2}(\lambda) = -(\lambda -k)[ \lambda^{2} - (a-c)\lambda -ac + \frac{1}{3} m^{2}].$
The roots of equation $ p_{2}(\lambda)= 0~$ are $~\lambda_{1}= k
,~\lambda_{2,3}= \frac{(a - c)\pm \sqrt{\Delta}}{2},$ where $ ~\Delta = (a+c)^{2} - \frac{4}{3} m^{2}.$\\
{\bf 1.} {\bf Case} $~k < 0~$ and $~\Delta < 0.~$  Then $~\lambda_{1}= k,~\lambda_{2,3}= \frac{(a-c)\pm i \sqrt{- \Delta}}{2}.~$ We have $~\Delta < 0~ $ if and only if $~m \in (-\infty, -\frac {|a+c|\sqrt {3}}{2} )\cup (\frac {|a+c|\sqrt {3}}{2} , \infty).$ \\
$ (i)~$ Assume that  $~k< 0~$ and $~a < c.~$ In this case  we have $ \lambda_{1} < 0 $ and $ Re (\lambda_{2,3})<
0.$ Since $|arg(\lambda_{i})| =\pi > \displaystyle\frac{q \pi}{2},
i=\overline{1,3}~$ for all $ q\in (0, 1)$, by Proposition 1(i), one obtains that $~e_{2}^{m}~$ is locally asymptotically stable.\\
$(ii)~$ Assume that $~k< 0~$ and $~a > c.~$  Then $~\lambda_{1} < 0 $ and $ Re(\lambda_{2,3})> 0.~$ In this case there exist two situations. Applying Proposition 1(i), we have:\\
{\bf (a)}$~$ $ e_{2}^{m}$ is locally asymptotically stable, for $~0 < q < q_{2},$ where $~ q_{2} = \displaystyle\frac{2}{\pi}
\arctan\displaystyle\frac{\sqrt{-\Delta}}{a-c}.~$  If $ q=q_{2},$ then $ e_{2}^{m} $ is stable.\\
 {\bf (b)}$~$ For $~ q_{2} < q < 1,~e_{2}^{m}$ is unstable.\\
{\bf 2.} {\bf Case} $~k < 0~$ and $~\Delta > 0.~$  Then $~\lambda_{1}= k,~\lambda_{2,3}= \frac{(a-c)\pm \sqrt{\Delta}}{2}\in {\bf R}.~$ We have $~\Delta > 0~ $ if and only if $~m \in (-\frac {|a+c|\sqrt {3}}{2}, \frac {|a+c|\sqrt {3}}{2}).$\\
$~(i)~$ Assume  $~k< 0.~$ Then $~\lambda_{1} <0.~$  We have $~\lambda_{2}<0~$ and $~\lambda_{3}<0~$  if and only if $~\lambda_{2} + \lambda_{3} <0 $ and  $~\lambda_{2}\cdot \lambda_{3} >0.~$ It follows that $~\lambda_{i} <0, i=\overline{1,3}~$ for $~a<c,~ ac > 0~$ and $~m\in (-\frac {|a+c|\sqrt {3}}{2},-\sqrt{3ac})\cup (\sqrt{3ac}, \frac {|a+c|\sqrt{3}}{2}).~$  According to Lemma 1(ii), $~ e_{2}^{m}$ is locally asymptotically stable for all $~q\in (0,1).$\\
 $~(ii)~$ Assume that  $~k< 0~$ and $~ ac < 0.~$  For all $~a,c \in {\bf R}^{\ast},~$  we have  $~\lambda_{2}>0~$ and $~\lambda_{3}>0~.$  Then, $~ e_{2}^{m}$ is unstable $~(\forall) m \in (-\sqrt{3ac}, \sqrt{3ac})~$ and $~q\in (0,1).$ \\
{\bf 3.} {\bf Case} $~k >0.~$ Let $~a, c\in {\bf R}^{\ast}.~$ Since $\lambda_{1} > 0,~ J_{k}(e_{2}^{m}, a, c)$ has at least
a positive eigenvalue. Then, by Lemma 1(i), $~e_{2}^{m} $ is unstable $ (\forall) m\in {\bf R}^{\ast}~$ and $~q\in (0,1).~$ \hfill$\Box$\\[-0.4cm]

{\bf 4.2.  Examples}

{\bf Example 1.} {\rm  Let be the special fractional model with control $ (15).$\\
$(1)~$ We select $~ a=-0.25, c  = 1 $ and  $ k= -0.25.~$ According to Theorem $ 3.1(i),~$ it follows that  $~ e_{0}= (0,0,0)~$ is
asymptotically stable for all $~ q\in (0, 1).$\\
$(2)~$ We consider $~ a=0.5, c =0.8, k= -0.75.~$  Applying Theorem $~3.1(ii),~$ it follows that $~ e_{0}=(0,0,0)~$ is unstable for all $~ q\in (0, 1).$}

{\bf Example 2.} {\rm  Let be the special fractional model with control  $ (15).$\\
$(1)~$ We select $~ a =-2,~c =1~$ and $~k=-0.8.~$  According to Theorem $~4.1(i),~ e_{2}^{m} = (0,m, 0)~$ is asymptotically stable $~(\forall) m\in (-\infty, -0,8861)\cup (0.8861, \infty)~ $ and  $~q\in (0,1).~$  In particular, taking $~m = 1~$ one obtains $~ \Delta=-\frac{1}{3}~$ and the eigenvalues  $~\lambda_{1} = - 0.8, ~ \lambda_{2,3}= -1.5\pm 0.8661 i.~$  Indeed,  $ e_{2} = (0,1,0)~$ is asymptotically stable for all $~q\in (0,1).~$\\
$(2)~$ We consider $~ a =2,~ c = 1~$ and $~k=- 2.~$ Taking $~m =\sqrt{7}\approx 2.6458, ~$ one obtains $~\Delta= - 0.(3)~$ and $ q_{2} = 0.(3).~$ According to Theorem $~4.1(ii)(a),~ e_{2} = (0,\sqrt{7}, 0)~$ is asymptotically stable $~(\forall) q\in (0, 0.(3))~.$
Indeed, $~\lambda_{1} = - 2<0, ~\lambda_{2,3}= 0.5\pm 0.2827 i~$ and $~m=\sqrt{7} \in (2.5981, \infty ).$ Then, $|arg(\lambda_{2,3})| =arctan \frac{0.2827}{0.5}= arctan~0.5654 > arctan~0.5233=\frac{\pi}{6} > q \frac{q \pi}{2},~(\forall) q\in (0, q_{2}),$ since $ \frac {\pi}{6} \approx 0.5233.$\\
$(3)~$ We take  $~ a =-1,~c =-0.5~$ and $~k=-0.4.~$  According to Theorem $~4.2(i),~ e_{2}^{m} = (0,m, 0)~$ is asymptotically stable $~(\forall) m\in (-1.2990, -1.2247)\cup (1.2247, 1.2990)~ $ and  $~q\in (0,1).~$  In particular, taking $~m =- 1.24~$ one obtains $~ \Delta= 0.1999~$ and  $~\lambda_{1} = - 0.4<0, ~\lambda_{2}= -0.0265<0~$ and $~\lambda_{3}= -0.4735<0.~$   Indeed,  $~ e_{2} = (0,-1.24, 0)~$ is asymptotically stable for all $~q\in (0,1).~$}\\

{\bf 5. Numerical integration of the special fractional system $ (15) $}\\[-0.3cm]

Let $f: {\bf R}\to {\bf R}$ be an integrable function. For $ q\in
(0, 1]$ we consider a function of
 $C^{1}-$class $~g^{q}: {\bf R}\times {\bf R} \to {\bf R},~(t,s)\to g_{t}^{q}(s)$
for $ t,s \in {\bf R}$ and $t_{0} \leq s \leq t.$

A {\it Riemann integral} of $f$ with respect to $ g_{t}^{q}$ is
defined by $~_{t_{0}} I_{t,g_{t}^{q}}f(t) =\int_{t_{0}}^{t} f(s)
g_{t}^{q}(s)ds.$

If the function $~g_{t}^{q}~$ is defined by\\[-0.2cm]
\begin{equation}
g_{t}^{q}(s) = \displaystyle\frac{1}{\Gamma(q)} (t-s)^{q-1}
e^{-\rho (t-s)},\label{(16)}
\end{equation}
where $\Gamma(q)$ is the Euler Gamma function and $\rho > 0,$ then
the Riemann integral $~_{t_{0}} I_{t,g_{t}^{q}}f(t)~$ becomes the
{\it fractional integral}, denoted with $_{t_{0}} I_{t}^{q}f(t)$,
  where:\\[-0.4cm]
\begin{equation}
_{t_{0}} I_{t}^{q}f(t) =
\displaystyle\frac{1}{\Gamma(q)}\int_{t_{0}}^{t}
(t-s)^{q-1}e^{-\rho (t-s)} f(s)ds.\label{(17)}
\end{equation}

{\bf Remark 2.} {\rm  For $\rho =0,$  the relation $(17)$ is the fractional
Riemann-Liouville integral (\cite{kilb}, Section 2.1). Also, the
fractional integral defined by $(17)$ is a special case of the
generalized fractional El-Nabulsi integral \cite{nabu}.}\hfill$\Box$

 If  $ \varphi:{\bf R}\times {\bf
R}^{n} \to {\bf R}^{n}$ is a deterministic function, we will
denote by:\\[-0.2cm]
\[
_{t_{0}} I_{t, g_{t}^{q}} \varphi(t,x(t)) =\int_{t_{0}}^{t}
\varphi(s,x(s)) g_{t}^{q}(s) ds.\\[-0.2cm]
\]

We call {\it Volterra integral equation with respect to $
g_{t}^{q},$} the functional Volterra type equation given by:\\[-0.4cm]
\begin{equation}
x(t) = x(t_{0}) +~ _{t_{0}}I_{t, g_{t}^{q}}
\varphi(t,x(t)).\label{(18)}
\end{equation}

The equation $(18)$ can be written formally in the following way:\\[-0.4cm]
\begin{equation}
dx(s) = \varphi(s,x(s))g_{t}^{q}(s)ds.\label{(19)}
\end{equation}

We apply the above considerations for the special fractional Chen-Lee system with control:\\[-0.4cm]
\begin{equation}
\left\{\begin{array}{lcl}
 D_{t}^{q} x^{i}(t) & = &
F_{i}(x^{1}(t), x^{2}(t), x^{3}(t)),~~~t\in
(t_{0}, \tau),~q \in (0,1)\\
x(t_{0}) &=& (x^{1}(t_{0}), x^{2}(t_{0}), x^{3}(t_{0}))
\end{array}\right.\label{(20)}\\[-0.2cm]
\end{equation}
where: $~F_{1}(x)=- x^{2} x^{3} + a x^{1},~ F_{2}(x)= x^{1} x^{3} + k  x^{2},~ F_{3}(x)= \frac{1}{3} x^{1}x^{2}-c x^{3}.$

Since the functions $ F_{i}(x(t)), i=\overline{1,3} $ are continuous, the initial value problem $(20)$ is equivalent to system of Volterra integral equations \cite{bora},
 which is given as follows:\\[-0.4cm]
\begin{equation}
x^{i}(t)~=~ x^{i}(t_{0})  + ~_{t_{0}}I_{t, g_{t}^{q}}
F_{i}(x^{1}(t), x^{2}(t), x^{3}(t)),
~~~~~i=\overline{1,3}.\label{(21)}
\end{equation}

The equations $(21)$ can be written in the following form:\\[-0.2cm]
\begin{equation}
 dx^{1}(s)= F_{1}(x(s)) g_{t}^{q}(s)ds
,~ dx^{2}(s)= F_{2}(x(s)) g_{t}^{q}(s)ds,~ dx^{3}(s)=
F_{3}(x(s))g_{t}^{q}(s)ds, \label{(22)}
\end{equation}
 where $g_{t}^{q}$ is given by $(16)$ and $
F_{i}(x(s)), i=\overline{1,3}$ are given in $(20).$

The system $(22)$ is called the {\it Volterra integral equations
associated to special fractional Chen-Lee system with control} $(15)$.

For the numerical integration of the system $(22)$ one can use
the first fractional Euler scheme \cite{difo, kilb}, which is expressed as follows:\\[-0.4cm]
\begin{equation}
x^{i}(j+1)=x^{i}(j)+ h F_{i}(x^{1}(j), x^{2}(j),
x^{3}(j))g_{t}^{q}(j),~~~ i=\overline{1,3},\label{(23)}
\end{equation}
where $ j=0,1,2,...,N,  h=\displaystyle\frac{T}{N}, T>0, N>0.~$

More precisely, if $t_{0}=0, \rho>0 $ and $ \varepsilon > 0, $ the numerical integration of the system $(23)$ is given by:\\[-0.4cm]
\begin{equation}
\left \{ \begin{array}{ll} x^{1}(j+1) &= x^{1}(j)+
h~\displaystyle\frac{1}{\Gamma(q)}(t-j)^{q-1}\varepsilon^{-\rho
(t-j)}(- x^{2}(j) x^{3}(j) + a x^{1}(j))\\[0.3cm]
x^{2}(j+1) &= x^{2}(j)+
h~\displaystyle\frac{1}{\Gamma(q)}(t-j)^{q-1}\varepsilon^{-\rho
(t-j)}( x^{1}(j) x^{3}(j) + k x^{2}(j))\\[0.3cm]
 x^{3}(j+1) &= x^{3}(j) +
h~\displaystyle\frac{1}{\Gamma(q)}(t-j)^{q-1}\varepsilon^{-\rho
(t-j)}(\frac{1}{3}
 x^{1}(j) x^{2}(j) - c x^{3}(j))\\[0.3cm]
x^{i}(0)&= x_{e}^{i}+\varepsilon,~~~i=\overline{1,3}.\\[0.3cm]
\end{array}\right. \label{(24)}
\end{equation}

{\bf Example 3.} {\rm  Let us we present the numerical simulation of solutions of the fractional system with control  $(15)$ which has considered in Example $~2.(1).~$
For this we apply the algorithm $ (24) $ and package Maple. Then, in $(24)$ we take: $~a = -2, c = 1, k =-0.8,~\rho =0.01,~h = 0.01, \varepsilon= 0.01,  N = 500, t = 502. $

The orbits $~Ox^{1}x^{2}x^{3}~$  for the solutions of the special fractional model $ (15) $ for
the equilibrium state $ e_{1} = (0, 1, 0)~ $ have the representations given in figures Fig. 1(a)(for $~q = 0.55 ~)$ and Fig. 1(b)(for $~ q = 1).$}\hfill$\Box$
\begin{center}
\begin{tabular}{ccc}
\includegraphics[width=5cm]{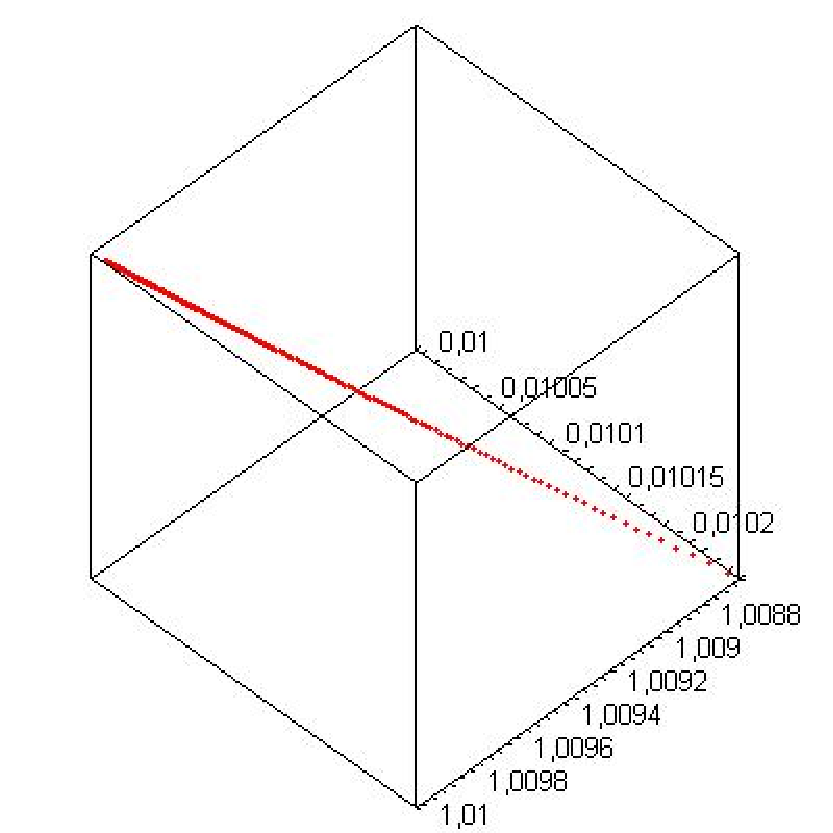}& &
\includegraphics[width=5cm]{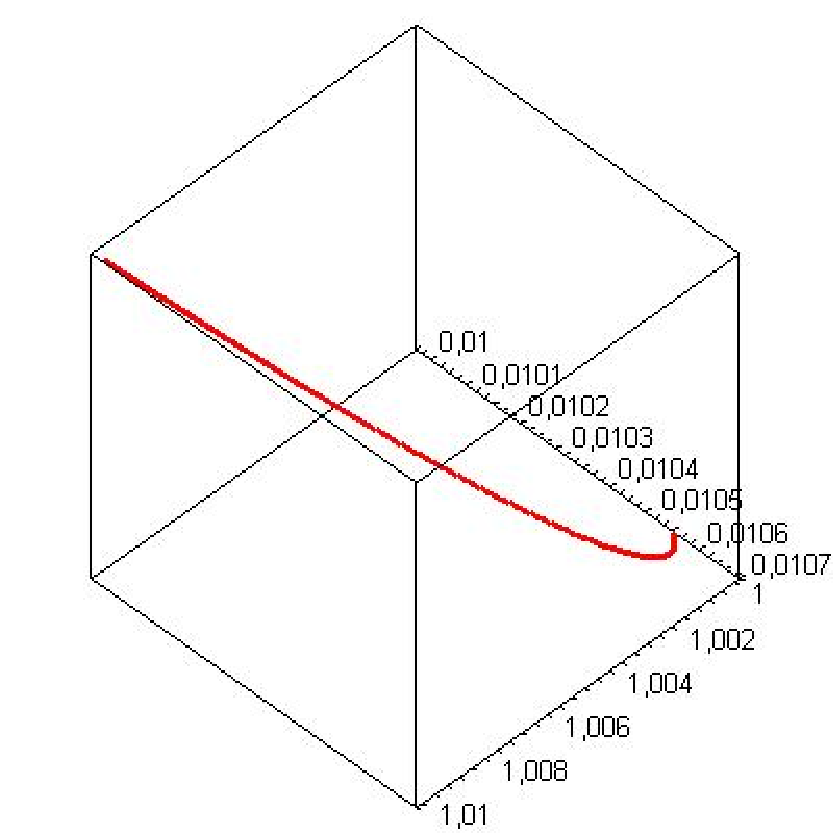}\\[-0.5cm]
\end{tabular}
\end{center}
\begin{center}
\begin{tabular}{lll}
{\bf Fig.1.(a)} $ (x^{1}(n), x^{2}(n), x^{3}(n)),~q=0.55$ & &
{\bf Fig.1.(b)} $ (x^{1}(n), x^{2}(n), x^{3}(n)),~q=1$\\[-0.5cm]
\end{tabular}
\end{center}

The numerical simulations confirm the validity of the theoretical analysis.\\[-0.3cm]

{\bf 6. Conclusions}\\[-0.4cm]

 The paper focuses on the special fractional Chen-Lee system for which the stability problem of its equilibrium states have been analyzed. Sufficient conditions on the parameters $~k, a $ and $ c $ in the special fractional model with control so that its equilibrium states are asymptotically stable have formulated and proved. Finally, we have described the numerical algorithm $ (24) $ in order to determine the approximate solution of system $ (15).$ The solutions of a special fractional model with control have been vizualized for a fixed vector $~v_{CL} $ (given in Example 2.(1)).

 By varying the values of parameters $~ k, a $ and $ c $ in system $ (15), $  it generates a series of chaotic and non-chaotic fractional dynamical systems. The studied fractional model with control can be applied in secure communications, complete synchronization and nonlinear circuit theory.\\[-0.7cm]

Author's adress\\[0.2cm]
\hspace*{0.7cm}West University of Timi\c soara. Seminarul de Geometrie \c si Topologie\\
\hspace*{0.7cm}  Teacher Training Department\\
\hspace*{0.7cm} Bd. V. P{\^a}rvan,no.4, 300223, Timi\c soara, Romania\\
\hspace*{0.7cm}E-mail: mihai.ivan@e-uvt.ro\\

\end{document}